\documentclass{article}
\author{Claire Levaillant\\clairelevaillant@yahoo.fr}
\title{Classification of the invariant subspaces of the
Lawrence--Krammer representation}
\date{}

\usepackage{amsmath, amssymb}
\usepackage{amsthm}

\newcommand{\n}{\nu}
\newcommand{\Q}{\mathbb{Q}}

\newcommand{\be}{\beta}
\newcommand{\al}{\alpha}
\newcommand{\xb}{x_{\beta}}
\newcommand{\e}{\epsilon}
\newcommand{\lra}{\longrightarrow}

\newcommand{\unsurr}{\frac{1}{r}}
\newcommand{\unsur}{\frac{1}}

\newcommand{\nts}{\negthickspace}
\newcommand{\ali}{\al_i}
\newcommand{\xali}{x_{\al_i}}

\newcommand{\V}{\mathcal{V}}
\newcommand{\IH}{\mathcal{H}}

\newcommand{\lb}{\lbrace}
\newcommand{\rb}{\rbrace}
\newcommand{\noin}{\noindent}

\newcommand{\ih}{\IH_{F,r^2}}

\newcommand{\W}{\mathcal{W}}
\newcommand{\di}{\text{dim}}

\newcommand{\da}{\downarrow}

\newcommand{\eg}{&=&}

\newcommand{\cil}{\frac{n(n-3)}{2}}
\newcommand{\dbw}{\frac{(n-1)(n-2)}{2}}
\newcommand{\chl}{\frac{n(n-1)}{2}}

\newcommand{\jca}{\unsur{r^{n-3}}}

\begin{document}
\maketitle
%\begin{center} \textit{Dedicated to Arjeh Cohen on his 60th
%birthday}\\\end{center}\vspace{0.2in}
\noin\textbf{Abstract.} The Lawrence--Krammer representation was
used in $2000$ to show the linearity of the braid group. The problem
had remained open for many years. The fact that the
Lawrence--Krammer representation of the braid group is reducible for
some complex values of its two parameters is now known, as well as
the complete description of these values. It is also known that when
the representation is reducible, the action on a proper invariant
subspace is an Iwahori-Hecke algebra action. In this paper, we prove
a theorem of classification for the invariant subspaces of the
Lawrence--Krammer space. We classify the invariant subspaces in
terms of Specht modules. We fully describe them in terms of
dimension and spanning vectors in the Lawrence--Krammer space.

\section{Introduction}\footnotetext{\\The author thanks Pierre-Albert Levaillant
for his support when part of this manuscript was written.} \noindent
The Lawrence--Krammer representation of degree $\frac{n(n-1)}{2}$
was used in $2000$ to show the linearity of the braid group on $n$
strands ($n\geq 3$). This result is due to Bigelow in \cite{BIG} and
independently to Krammer in \cite{KR}. Their proofs are very
different. Krammer's proof is algebraic while Bigelow's proof is
topological. Linearity of a group means that there exists a faithful
linear representation of this group. The Lawrence--Krammer
representation first appears in a work of Lawrence in \cite{RUT}. It
is thus called the Lawrence--Krammer representation. To show that
the braid goup on $n$ strands is a linear group, the $n$-dimensional
Burau representation was long a candidate. However, the Burau
representation is unfaithful for $n\geq 5$, see \cite{MOO},
\cite{LP}, \cite{BI}. Krammer's representation and proof of
linearity was generalized by Cohen-Wales in \cite{CW} for Artin
groups of finite type. The same result of linearity is proven
independently by Digne in \cite{DI}.
\\The Lawrence--Krammer representation based on two parameters $t$
and $q$ was known by several authors (Cohen--Gijsbers--Wales in
\cite{CGW}, Marin in \cite{MAR}, Zinno in \cite{Z}) to be
generically irreducible. It is shown in \cite{LEV} with some
restrictions on the parameter $q$ that when the two parameters are
specialized to some nonzero complex numbers, the representation
becomes reducible. The complete list of the nonzero complex
parameters for which the Lawrence--Krammer representation is
reducible is given in \cite{LEV}. In the same paper, it is shown
that when the Lawrence--Krammer representation is reducible, the
action on a proper invariant subspace of the Lawrence--Krammer space
is an Iwahori-Hecke algebra action. In this paper, we give the
complete classification of the proper invariant subspaces of the
Lawrence--Krammer space in terms of Specht module. Our result is as follows. It was originally stated as a conjecture in the author's thesis in \cite{THE}.
This result also follows independently from \cite{RUI} where Rui and Si use the notions of cellular algebras and Jucys-Murphy bases, which we don't use here.\\
\newtheorem*{thm}{Main Theorem}
\begin{thm}
Let $n$ be an integer with $n\geq 3$. Assume that $q$ is not a
$k$-th root of unity for every integer $k$ with $1\leq k\leq n$. There are two cases:\\
$1)$ Assume that $q^n\neq -1$ when $t=-1$: \\
if the Lawrence--Krammer representation is reducible, its unique
proper invariant subspace is isomorphic to one of the Specht modules
$$S^{(n)},\,S^{(n-1,1)},\,S^{(n-2,2)},\,S^{(n-2,1,1)},$$
which respectively arise if and only if
$$t=\frac{1}{q^n},\,t\in\bigg\lbrace\frac{1}{\sqrt{q}^n},-\frac{1}{\sqrt{q}^n}\bigg\rbrace,\,t=\frac{1}{q},\,t=-1$$
$2)$ If $t=-1$ and $q^n=-1$, there are exactly three proper
invariant subspaces of the Lawrence--Krammer space and they are
respectively isomorphic to $S^{(n)}$, $S^{(n-2,1,1)}$ and
$S^{(n)}\oplus S^{(n-2,1,1)}$\\
\end{thm}

\noin In the next sections, we introduce the Lawrence--Krammer
representation and prove the Main Theorem. Further, we fully
describe all the proper invariant subspaces by providing their
dimensions and some spanning vectors in the Lawrence--Krammer space.

\section{The Lawrence--Krammer representation}
The Lawrence--Krammer space $\V^{(n)}$, abbreviated L-K space, is
the vector space of dimension $\frac{n(n-1)}{2}$ over the field
$\mathbb{Q}(t,q)$ with spanning vectors the $\xb$'s, indexed by the
$\frac{n(n-1)}{2}$ positive roots $\be$'s of a root system of type
$A_{n-1}$. In what follows, we will denote the simple roots by
$\al_1,\dots,\al_{n-1}$. The positive roots are the sums
$\al_i+\dots+\al_j$ (with $i\leq j$) of simple roots. If
$\be=\al_i+\dots+\al_j$, we will denote the basis vector $\xb$ by
$w_{i,j+1}$. The height $ht(\be)$ of a positive root $\be$ is the
sum of its coefficients with respect to the simple roots. These
coefficients are either zeros or ones. The support $Supp(\be)$ of a
positive root $\be$ is the set of $k\in\lb 1,\dots,n\rb$ such that
the coefficient of $\al_k$ in $\be$ is nonzero. We will denote the
set of positive roots by $\phi^{+}$. Finally, if $(m_{ij})_{1\leq
i\leq j\leq n-1}$ denotes the Coxeter matrix of type $A_{n-1}$, the
inner product between two simple roots $\al_i$ and $\al_j$ is given
by
$$(\al_i|\al_j)=-cos\bigg(\frac{\pi}{m_{ij}}\bigg)$$ \noin The Lawrence--Krammer
representation can be constructed via BMW algebras. These algebras
are algebras named after Birman, Murakami and Wenzl. They were
introduced in \cite{BW} in order to study the linearity of the braid
groups and independently by Murakami in \cite{MUR}. They feature in
many areas including statistical mechanics, knot theory and quantum
groups theory. The BMW algebra $B(A_{n-1})$ or simply $B$ of type
$A_{n-1}$ with parameters $l$ and $m$ as defined in \cite{CGW} is
the algebra over the field $\mathbb{Q}(l,m)$ with $(n-1)$ generators
$g_1,\dots,g_{n-1}$, those of the braid group, and other elements
$e_1,\dots,e_{n-1}$ that are related to the $g_i$'s by
$$m\,e_i=l(g_i^2+m\,g_i-1)$$

\noin The other defining relations that relate the elements $e_i$'s
and $g_i$'s of the algebra are the following.
\begin{eqnarray*}
g_ie_i&=&l^{-1}e_i\qquad\text{for all $i$}\\
e_ig_je_i&=&l\,e_i\qquad\;\;\,\text{when $|i-j|=1$}
\end{eqnarray*}

\noin The BMW algebra $B$ modulo the two-sided ideal $I_1=Be_1B$ is
the Hecke algebra with generators $g_1,\dots,g_{n-1}$ and relations
the braid relations and the relations
$$g_i^2+m\,g_i=1\qquad\text{for all $i$}$$\noin When $m$ is a given
nonzero complex number, we let $r$ and $-\unsurr$ be the two nonzero
complex roots of the quadratic $X^2+m\,X-1=0$. So, the nonzero
complex numbers $m$ and $r$ are related by $m=\unsurr-r$.\\ Up to a
rescaling of the generators, the algebra $B/I_1$ is the
Iwahori-Hecke algebra $\ih(n)$ of the symmetric group $Sym(n)$ with
parameter $r^2$ over the field $F=\Q(l,r)$, as defined in
\cite{MAT}.

\noin In \cite{THE}, the author uses the isomorphism between the BMW
algebra and the tangle algebra of Morton-Traczyk (see \cite{MW}) to
construct a representation $\n^{(n)}$ of degree $\chl$ of the BMW
algebra of type $A_{n-1}$ inside the Lawrence--Krammer space
$\V^{(n)}$ over $F$. She shows that as a representation of the braid
group on $n$ strands and up to some change of parameters and some
rescaling of the generators, this representation is equivalent to
the Lawrence--Krammer representation. The change of parameters is
given by $lt=r^3$ and $q=\unsur{r^2}$. The representation $\n^{(n)}$
is defined on the elements $g_i$'s and $e_i$'s of the algebra by
$$\n^{(n)}:\left.\begin{array}{l}
B  \lra  End_{F}(\V^{(n)})\\
g_i  \longmapsto  \n_i\\
e_i  \longmapsto
\n^{(n)}(e_i)=\frac{l}{m}(\n_i^2+m\,\n_i-id_{\V^{(n)}})
\end{array}\right.$$
where $\n_i$ is the endomorphism defined on the basis vectors
$\xb$'s by $$\nts\nts\nts\nts\nts\nts\n_i(\xb)=\begin{cases}
r\,\xb&\text{if $\left.\begin{array}{l} Supp(\be)\cap\lb i-1,i,i+1\rb=\emptyset\\
\text{or}\;\,\lb i-1,i,i+1\rb\subseteq Supp(\be)\end{array}\right.$}\\
\unsur{l}\,\xb&\text{if $\;\;\be=\ali$}\\
x_{\be+\al_i}&\text{if $\be=\al_s+\dots+\al_{i-1}$ with $s\leq i-1$ $(a)$}\\
x_{\be+\al_i}+m\,r^{ht(\be)-1}\,\xali-m\,\xb&\text{if $\be=\al_{i+1}+\dots+\al_k$ with $k\geq i+1$ $(b)$}\\
x_{\be-\al_i}+\frac{m}{l\,r^{ht(\be)-2}}\,\xali-m\,\xb&\text{if $\be=\al_s+\dots+\ali\;\;\;\,$ with $s\leq i-1$ $\,(c)$}\\
x_{\be-\al_i}&\text{if $\be=\ali+\dots+\al_k\;\;\;\,$ with $k\geq
i+1$ $(d)$}
\end{cases}$$
Note that $(a)$ and $(b)$ are the two different ways in which the
inner product $(\ali|\be)$ can be $-\unsur{2}$ and $(c)$ and $(d)$
are the two different ways in which the inner product $(\ali|\be)$
can be
$\unsur{2}$. \\
We now show that the representation $\n^{(n)}$ is equivalent to the
Lawrence--Krammer representation of the BMW algebra defined by
Cohen--Gijsbers--Wales. In \cite{CGW}, the authors define $I_2$ as
the two-sided ideal of $B$ generated by all the products $e_ie_j$
with $|i-j|>1$. For each irreducible representation $\theta$ of the
Hecke algebra of type $A_{n-3}$, they construct a representation of
$B/I_2$ of degree $|\phi^{+}|deg(\theta)$ and they show that these
are all the irreducible representations of $I_1/I_2$. The
Lawrence--Krammer representation of the BMW algebra is obtained for
one the two inequivalent representations of degree one of the Hecke
algebra of type $A_{n-3}$. We next show that $\n^{(n)}$ is
irreducible and factors through $I_2$. The proof of the following
result can be found in \cite{THE}.
\newtheorem*{Theor}{Theorem}
\begin{Theor} Assume $r^{2k}\neq 1$ for every integer $k$ with
$1\leq k\leq n$. \\When $n\geq 4$, $\n^{(n)}$ is irreducible except
when $l\in\lb
r,-r^3,\unsur{r^{2n-3}}, -\unsur{r^{n-3}}, \unsur{r^{n-3}}\rb$ when it is reducible.\\
$\n^{(3)}$ is irreducible except when $l\in\lb -r^3,\unsur{r^3}, -1,
1\rb$ when it is reducible.\end{Theor} \noin This result shows in
particular that $\n^{(n)}$ is generically irreducible for every
integer $n$ with $n\geq 3$. Moreover, straightforward computations
show that for every nodes $i$ and $j$ with $|i-j|>1$, we have
$\n^{(n)}(e_ie_j)=0$. Thus, $\n^{(n)}$ is an irreducible
representation of $I_1/I_2$ of degree $|\phi^{+}|$. It must then be
equivalent to the Lawrence--Krammer representation of the BMW
algebra of \cite{CGW}. Our $r$ is the $\unsurr$ of \cite{CGW}.\\
We note that the restriction on $r$ that $r^{2k}\neq 1$ for every
integer $k$ with $1\leq k\leq n$ is equivalent to assuming that
$\ih(n)$ is semisimple (see Corollary $3.44$ page $48$ of
\cite{MAT}).

\section{The invariant subspaces of the L-K space}

When the representation $\n^{(n)}$ is reducible, the action on a
proper invariant subspace of $\V^{(n)}$ is an Iwahori-Hecke algebra
action: this is Proposition $1$ of \cite{CLW}. The following two
theorems stated here for the Iwahori-Hecke algebra of the symmetric
group $Sym(n)$ instead of the symmetric group $Sym(n)$ are due to
James in \cite{JAM}. In characteristic zero, when the Iwahori-Hecke
algebra of the symmetric group is semisimple, they remain true for
the Iwahori-Hecke algebra of the symmetric group (see \cite{MAT}).
\newtheorem{theo}{Theorem}
\begin{theo}
Let $n$ be an integer with $n\geq 7$ and assume that $\ih(n)$ is
semisimple. \\
Then, every irreducible $\ih(n)$-module is either isomorphic to one
of the Specht modules $S^{(n)}$, $S^{(1^n)}$, $S^{(n-1,1)}$,
$S^{(2,1^{n-2})}$ or has dimension greater than $(n-1)$.
\end{theo}
\begin{theo} Let $n$ be an integer with $n\geq 9$ and assume that $\ih(n)$ is semisimple.\\
Then, every irreducible $\ih(n)$-module is either isomorphic to one
of the Specht modules $S^{(n)}$, $S^{(n-1,1)}$, $S^{(n-2,2)}$,
$S^{(n-2,1,1)}$ or their conjugates, or has dimension greater than
$\frac{(n-1)(n-2)}{2}$.
\end{theo}

\noin Theorem $1$ fails for $n=4$ as $S^{(2,2)}$ has dimension $2$
and for
$n=6$ as $S^{(3,3)}$ and $S^{(2,2,2)}$ both have dimensions $5$. \\
Theorem $2$ fails for $n=7$ as $S^{(4,3)}$ has dimension $14$ and
for $n=8$ as $S^{(4,4)}$ and $S^{(2,2,2,2)}$ have dimensions $14$.
\\
A consequence of Theorem $1$ and Theorem $2$ is the following
corollary.

\newtheorem{Cor}{Corollary}
\begin{Cor}
Let $n$ be an integer with $n\geq 5$ and $n\neq 8$ and assume that $\ih(n)$ is semisimple. \\
Then, the irreducible $\ih(n)$-modules have dimension $1$, $n-1$,
$\cil$, $\dbw$ or dimension greater than $\dbw$.\\
Assume that $\ih(8)$ is semisimple. Then, the irreducible
$\ih(8)$-modules have dimension $1$, $7$, $14$, $20$, $21$ or
dimension greater than $21$.
\end{Cor}

\noin We now recall some results of \cite{CLW} about the existence
of a one-dimensional invariant subspace of the L-K space and of an
irreducible $(n-1)$-dimensional invariant subspace of the L-K space
for some values of the parameters $l$ and $r$.

\begin{theo}
Let $n$ be an integer with $n\geq 3$ and assume $\ih(n)$ is semisimple.\\
\indent Suppose $n\geq 4$. There exists a one-dimensional invariant
subspace of $\V^{(n)}$ if and only if $l=\unsur{r^{2n-3}}$. If so,
it is spanned by $$\sum_{1\leq s<t\leq n}r^{s+t}\,w_{st}$$
\indent\textit{Case $n=3$.} There exists a one-dimensional invariant
subspace of $\V^{(3)}$ if and only if $l=\unsur{r^3}$ or
$l=-r^3$.\vspace{0.05in} \\Moreover, if $r^6\neq -1$, it is unique
and \vspace{-0.099in}
$$\begin{array}{cccccc} \text{when}&l&=&\unsur{r^3},&\text{it is
spanned
by}&w_{12}+r\,w_{13}+r^2\,w_{23}\\
\text{when}&l&=&-r^3,&\text{it is spanned
by}&w_{12}-\unsurr\,w_{13}+\unsur{r^2}\,w_{23}
\end{array}\vspace{-0.11in}$$
\noin If $r^6=-1$, there are exactly two one-dimensional invariant
subspaces of $\V^{(3)}$ and they are respectively spanned by the
vectors above.
\end{theo}
\noin \textbf{Proof.} This is Theorem $4$ of \cite{CLW}. It is in
particular shown along the proof that except when $n=3$, the Specht
module $S^{(n)}$ occurs in the L-K space $\V^{(n)}$ for
$l=\unsur{r^{2n-3}}$, while the conjugate Specht module $S^{(1^n)}$
cannot occur in the L-K space. When $n=3$, $S^{(3)}$ occurs for the
value $l=\unsur{r^3}$ and $S^{(3-2,1,1)}$ occurs for the value
$l=-r^3$.

\begin{theo} Let $n$ be an integer with $n\geq 3$ and $n\neq 4$. Assume
$\ih(n)$ is semisimple. Then, there exists an irreducible
$(n-1)$-dimensional invariant subspace of $\V^{(n)}$ if and only if
$l=\unsur{r^{n-3}}$ or
$l=-\unsur{r^{n-3}}$. \\\\
If so, it is spanned by the $v_i^{(n)}$'s, $1\leq i\leq n-1$, where
$v_i^{(n)}$ is defined by the formula:
\begin{equation*}\begin{split}v_i^{(n)}=\Big(\unsur{r}\,-\,\unsur{l}\Big)w_{i,i+1}&+\sum_{s=i+2}^{n}r^{s-i-2}(w_{i,s}\,-\,\unsurr\;w_{i+1,s})\\
&+\;\e_l\,\sum_{t=1}^{i-1}r^{n-i-2+t}(w_{t,i}\,-\,\unsurr\;w_{t,i+1})\end{split}\end{equation*}
\begin{center} with \hspace{.1in}
$\begin{cases}\e_{\unsur{r^{n-3}}}\;\;=\,1\\\e_{-\unsur{r^{n-3}}}=-1\end{cases}$
\end{center}
\noin\textit{Case $n=4$.} Assume $\ih(4)$ is semisimple. Then, there
exists an irreducible $3$-dimensional invariant subspace of
$\V^{(4)}$ if and only if
$l\in\lb\unsurr,-\unsurr,-r^3\rb$.\\\\
If $l\in\lb -\unsurr,\unsurr\rb$, it is spanned by $v_1^{(4)}$, $v_2^{(4)}$, $v_3^{(4)}$.\\\\
If $l=-r^3$, it is spanned by the vectors:
$$\left\lb\begin{array}{cccccccccccccc} u_1\eg r\,w_{23}+
w_{13}+(\unsurr+\unsur{r^3})w_{34}-
w_{24}- \unsurr\,w_{14}\\
u_2\eg -r\,w_{12}-r^2\,w_{13}-
\unsurr\,w_{34}-\unsur{r^2}\,w_{24}+(r+\unsurr)\,w_{14}\\
u_3\eg (r+r^3)\,w_{12}+\unsurr\,w_{23}- w_{13}+w_{24}-
r\,w_{14}\\
\end{array}\right.$$
\end{theo}

\noin\textbf{Proof.} This is Theorem $5$ of \cite{CLW}. In
particular, it is shown that for $n\geq 3$ and $n\neq 4$,
$S^{(n-1,1)}$ occurs in the L-K space $\V^{(n)}$ for
$l=\unsur{r^{n-3}}$ and for $l=-\unsur{r^{n-3}}$, while the
conjugate Specht module $S^{(2, 1^{n-2})}$ cannot occur in the L-K
space. When $n=4$, $S^{(3,1)}$ occurs in $\V^{(4)}$ when
$l\in\lb\unsurr,-\unsurr\rb$ and $S^{(4-2,1,1)}$ occurs in
$\V^{(4)}$ for $l=-r^3$. In the same proof, it is also shown that
the Specht modules $S^{(3,3)}$ and its conjugate $S^{(2,2,2)}$ both
of dimension $5$ cannot occur inside $\V^{(6)}$.
\newtheorem{Remark}{Remark}
\begin{Remark} For $n\geq 5$, Theorem $3$, Theorem $4$ and Corollary $1$
imply that if $\n^{(n)}$ is reducible and if
$l\not\in\lb\unsur{r^{n-3}},-\unsur{r^{n-3}},\unsur{r^{2n-3}}\rb$,
then an irreducible invariant subspace of $\V^{(n)}$ must have
dimension greater than or equal to $\cil$ when $n\neq 8$ and greater
than or equal to $14$ when $n=8$.\end{Remark} \noin In fact, we have
the following theorem.
\begin{theo}
Let $n$ be an integer with $n\geq 4$ and assume that $\ih(n)$ is
semisimple. Then, there exists an irreducible $\cil$-dimensional
invariant subspace of $\V^{(n)}$ if and only if $l=r$. If so it is
unique.
\end{theo}
\noin \textbf{Proof.} The idea is to consider the $B$-module
$K(n)=\cap_{1\leq i<j\leq n} Ker\,\n^{(n)}(C_{ij})$, where
$C_{i,i+1}=e_i$ and $C_{ij}=g_{j-1}^{-1}\dots
g_{i+1}^{-1}e_ig_{i+1}\dots g_{j-1}$ for all $j\geq i+2$. We will
denote its dimension by $k(n)$.
\begin{Remark} Since the $e_i$'s act trivially on any proper invariant subspace $\W$ of $\V^{(n)}$,
such a space $\W$ must be contained in $K(n)$.\end{Remark}
%As any proper invariant subspace $\W$ of $\V^{(n)}$ is
%5annihilated by the $e_i$'s, it must also be annihilated by the $g_k$
%conjugates of the $e_i$'s. Hence such a space $\W$ must be contained
%in $K(n)$. \end{Remark}
\noindent When $l=r$, we show that $K(n)$ is irreducible. %This is
%part of Proposition $6$, page $99$ of \cite{AR}.
For $n=4$, this result is part of Proposition $3$ of \cite{CLW}.
When $n\geq 5$, Proposition $4$ of \cite{CLW} shows that $K(n)$ is
nonzero. To show that $K(n)$ is irreducible, the idea is to use the
fact that $K(n)$ is an $\ih(n)$-module. Suppose first $n\neq 8$. If
$K(n)$ is reducible, by semisimplicity of $\ih(n)$, the
$\ih(n)$-module $K(n)$ decomposes as a direct sum $K_1(n)\oplus
K_2(n)$ of $\ih(n)$-modules with $K_1(n)$ irreducible. Since we have
$\chl-\dbw=n-1$ and $\cil=\dbw-1$, Corollary $1$ implies that one of
the two modules $K_1(n)$ or $K_2(n)$ must have dimension less than
or equal to $(n-1)$. We now recall that a necessary and sufficient
condition on $r$ so that $\ih(n)$ is semisimple is that $r^{2k}\neq
1$ for every integer $k$ with $1\leq k\leq n$. In particular, when
$l=r$, we have $l\not\in\lb\jca, -\jca,\unsur{r^{2n-3}}\rb$. Then,
by Remark $1$, it is impossible to have $K_1(n)$ or $K_2(n)$ of
dimension less than or equal to $(n-1)$. Thus, $K(n)$ is irreducible
and $k(n)\geq\cil$ by Remark $1$. When $n=8$ the proof is the same
but needs to be slightly adapted. We obtain that $K(8)$ is
irreducible and $k(8)\geq 14$. A consequence of Remark $2$ and of
the irreducibility of $K(n)$ is that $K(n)$ is the unique proper
invariant subspace of $\V^{(n)}$. For $n=4$, the result of Theorem
$5$ is Proposition $3$ of \cite{CLW}. When $n\geq 5$ and $l=r$, it
is shown in \cite{LEV} that $k(n)\leq\cil$ (see proof of Theorem
$3.3$. For a detailed proof, see pages $112-116$ of \cite{AR}).
Then, when $n\geq 5$ and $n\neq 8$, $k(n)=\cil$. When $n=8$, we have
$k(8)\in\lb 14,20\rb$. But if $k(8)=14$, then we observe that
$k(8)=k(7)$. We will use the following lemma to get a contradiction.
In this lemma, the Iwahori-Hecke algebra $\ih(n)$ is still assumed
to be semisimple.
\newtheorem{Lemma}{Lemma}
\begin{Lemma}
Suppose $l=r$. Then $K(n-1)\subseteq K(n)$ for all $n\geq 5$.
\end{Lemma}
\noin \textbf{Proof of the lemma.} Let $n\geq 5$. %If $K(n)\cap\V^{(n-1)}=\lb
%0\rb$, then $K(n)\oplus\V^{(n-1)}\subseteq\V^{(n)}$, which implies
%$k(n)\leq n-1$. Suppose first $n\neq 8$. Then $k(n)\geq\cil$ and for
%$n\geq 5$, we have $\cil>n-1$, hence a contradiction. Thus,
%$K(n)\cap\V^{(n-1)}\neq \lb 0\rb$. Moreover,
The vector space $K(n)\cap\V^{(n-1)}$ is not the whole space
$\V^{(n-1)}$ (for a proof, see the arguments of the proof of
Proposition $1$ of \cite{CLW}). Then by Remark $2$, we have
$K(n)\cap\V^{(n-1)}\subseteq K(n-1)$. Moreover, by Proposition $5$,
Chapter $8$ of \cite{THE}, we have $K(n)\cap\V^{(n-1)}\neq\lb 0\rb$.
Hence by irreducibility of $K(n-1)$, we must have
$K(n)\cap\V^{(n-1)}=K(n-1)$, which implies in
particular $K(n-1)\subseteq K(n)$.\\

\noin Let's go back to the proof of Theorem $5$. By the lemma, we
get $K(8)=K(7)$. By Proposition $5$, Chapter $8$ of \cite{THE}, the
element $r^2\,w_{12}-r\,w_{13}+w_{34}-r\,w_{24}$ belongs to $K(8)$.
We act with $\n_7\dots\n_4$ to see that the element
$r^4(r^2\,w_{12}-r\,w_{13})+w_{38}-r\,w_{28}$ also belongs to
$K(8)$. However, this element is not in $K(7)$, so we get a
contradiction. Hence it is impossible to have $k(8)=14$ and so
$k(8)=20$. Thus, for all $n\geq 4$, we have shown that when $l=r$,
the $B(A_{n-1})$-module $K(n)$ is the
unique proper invariant subspace of $\V^{(n)}$ and it has dimension $\cil$.\\
Conversely, it is shown in \cite{LEV} (see proof of Theorem $3.3$
and forthcoming \cite{dim}) that if there exists an irreducible
$\cil$-dimensional invariant subspace of the L-K space, then $l=r$.
\\\\We now describe the irreducible $\cil$-dimensional invariant
subspace $K(n)$ of $\V^{(n)}$ when $l=r$.

\begin{theo}
\hfill\\

\noin Assume $l=r$. Let $n$ be an integer with $n\geq 4$ and assume
that $\ih(n)$ is semisimple. \begin{itemize}\item When $n=4$, the
unique invariant subspace $K(4)$ of $\V^{(4)}$ is spanned by the two
linearly independent vectors:
$$\begin{array}{ccc}
w_1^{(4)}&=&(w_{14}-\unsurr\,w_{24})+(w_{23}-r\,w_{13})\\
w_2^{(4)}&=&(w_{24}-\unsurr\,w_{34})+(w_{13}-r\,w_{12})
\end{array}$$
\item When $n\geq 5$, the unique invariant subspace $K(n)$ of $\V^{(n)}$ is
built inductively as a direct sum of the unique invariant subspace
$K(n-1)$ of $\V^{(n-1)}$ and of an $(n-2)$-dimensional vector space
spanned by the vectors:
$$\begin{array}{cccccc}
w_1^{(n)}&=&w_{1,n}-\unsurr\,w_{2,n}&+&r^{n-4}\,(w_{23}-r\,w_{13})&\\
w_k^{(n)}&=&w_{k,n}-\unsurr\,w_{k+1,n}&+&r^{n-4}\,(w_{1,k+1}-r\,w_{1,k}),&
2\leq k\leq n-2
\end{array}$$
\end{itemize}
\end{theo}
\noin \textbf{Proof.} When $n=4$, see Proposition $3$ of \cite{CLW}.
When $n\geq 5$, we have seen that $K(n-1)\subseteq K(n)$. Hence, it
suffices to check that the $(n-2)$ linearly independent vectors of
the theorem belong to $K(n)$. This is achieved in \cite{THE},
Chapter $10$.\\\\
Let's now study the case of reducibility $l=-r^3$. This case
requests more attention. Indeed, when $r^{2n}=-1$, we have
$l=-r^3=\unsur{r^{2n-3}}$. In that case, $K(n)$ is no longer
irreducible. In fact we have the following result.

\begin{theo}
Let $n$ be an integer with $n\geq 5$. Assume $\ih(n)$ is
semisimple.\begin{enumerate}\item If $l=-r^3$ and $r^{2n}\neq -1$,
then $K(n)$ is irreducible and $k(n)=\dbw$. In particular, $K(n)$ is
the unique proper invariant subspace of $\V^{(n)}$.
\item If $l=-r^3$ and $r^{2n}=-1$, then $K(n)$ is reducible and
$k(n)=1+\dbw$. Moreover, $K(n)$ is a direct sum of the unique
one-dimensional invariant subspace of $\V^{(n)}$ and of the unique
irreducible $\dbw$-dimensional invariant subspace of
$\V^{(n)}$.\end{enumerate}
\end{theo}

\noin\textbf{Proof.} When $r^{2n}\neq -1$, the proof of
irreducibility of $K(n)$ is the same as in the case $l=r$. Moreover,
by Lemma $10$, Chapter $9$ of \cite{THE}, we know that
$k(n)\leq\dbw$. If $n\neq 8$, we hence get $k(n)=\dbw$. When $n=8$,
if $k(8)=14$ it comes $7\leq\text{dim}(K(8)\cap\V^{(7)})\leq 14$.
But since $-r^3\not\in\lb r,-\unsur{r^4},\unsur{r^4}\rb$, this is impossible. Hence the case $n=8$ is not an exception and $k(8)=21$.\\
When $r^{2n}=-1$, we have $l=-r^3=\unsur{r^{2n-3}}$. Then, there
exists a one-dimensional invariant subspace of $\V^{(n)}$ by Theorem
$3$. We have $0\subset K(n)\cap\V^{(n-1)}\subseteq K(n-1)$, where
the first inclusion holds by Proposition $5$ of \cite{THE}, Chapter
$8$. Moreover, since $r^{2n}=-1$, we have $r^{2(n-1)}\neq -1$. Thus
by the first point, $K(n-1)$ is irreducible. Hence it comes
$K(n)\cap\V^{(n-1)}=K(n-1)$. If $K(n)$ were one-dimensional, so
would be $K(n)\cap\V^{(n-1)}$. Then $K(n-1)$ would also be
one-dimensional. This would force $l=\unsur{r^{2n-5}}$, which is
impossible. Hence the one-dimensional invariant subspace of
$\V^{(n)}$ has a summand $S$ in $K(n)$. In particular $K(n)$ is
reducible. Moreover, except possibly when $n=8$, the uniqueness part
in Theorem $3$ and Theorems $4$ and $5$ allow to claim that the
summand $S$ has dimension greater than or equal to $\dbw$. As for
$n=8$, if $\text{dim}(S)=14$, then it comes $k(8)=15$. By arguments
already discussed before, we have $K(7)\subseteq K(8)$. Since by the
first point we have $k(7)=15$, it follows that $K(7)=K(8)$. When
$l=-r^3$, the vector $$-r\,w_{23}-\unsurr\,w_{34}+w_{24}$$ belongs
to $K(8)$ by Proposition $5$ of Chapter $8$ of \cite{THE}. By acting
with $\n_7\n_6\n_5\n_4$ on this vector, we see that
$$-r^5\,w_{23}-\unsurr\,w_{38}+w_{28}$$ also belongs to $K(8)$. The
latter vector is not in $K(7)$, hence a contradiction.
%But by the
%same arguments as before, since $r^{14}\neq -1$, we have
%$K(7)\subset K(8)$. Then $k(7)<15$. However by the first point,
%$k(7)=15$, hence a contradiction.
So in any case, we have $k(n)\geq 1+\dbw$. Further, from
$K(n)\cap\V^{(n-1)}=K(n-1)$, we derive $k(n)\leq k(n-1)+(n-1)$.
Replacing $k(n-1)=\frac{(n-2)(n-3)}{2}$, we get $k(n)\leq 1+\dbw$.
Gathering both inequalities now yields $k(n)=1+\dbw$. From this
equality on the dimensions, we deduce the existence of an
irreducible $\dbw$-dimensional invariant subspace of $\V^{(n)}$. It
remains to show that it is unique. Let $\W$ be an irreducible
$\dbw$-dimensional invariant subspace of $\V^{(n)}$ such that $K(n)$
is a direct sum of $\W$ and of the one-dimensional invariant
subspace of $\V^{(n)}$. Since for $n\geq 5$, we have $\dbw>n-1$, it
follows that the intersection $\W\cap\V^{(n-1)}$ is nontrivial.
Then, by irreducibility of $K(n-1)$, we get
$\W\cap\V^{(n-1)}=K(n-1)$. In particular, $K(n-1)\subseteq \W$. Let
$S$ be an $\ih(n-1)$ summand of $K(n-1)$ in $\W$. Since
$k(n-1)=\frac{(n-2)(n-3)}{2}$, this summand must be
$(n-2)$-dimensional. To conclude, it will suffice to prove the
following two lemmas.
\begin{Lemma}
Let $n$ an integer with $n\geq 5$. Assume that $\ih(n)$ is
semisimple. Suppose $l=-r^3$ and $r^{2n}=-1$. In $K(n)$ there exists
a unique one-dimensional $\ih(n-1)$-module, namely the unique
one-dimensional invariant subspace of $\V^{(n)}$.
\end{Lemma}
\noin\textbf{Proof.} This is an adaptation of Lemma $13$ page $139$
of \cite{AR}, where the assumption $r^{2(n-1)}=-1$ on $r$ has been
replaced with the assumption $r^{2n}=-1$. The scalar $\mu$ of
\cite{AR} must then take the value $1$ instead of the value $0$.
Thus, if such a space exists, it must be spanned by $$\sum_{1\leq
i<j\leq n}r^{i+j}\,w_{ij}$$ This ends the proof of the lemma. A
consequence of this lemma is that $S$ is an irreducible
$\ih(n-1)$-module of dimension $(n-2)$.
\begin{Lemma}
Let $n$ be an integer with $n\geq 5$. Assume that $\ih(n)$ is
semisimple. Suppose $l=-r^3$ and $r^{2n}=-1$. In $K(n)$ there exists
a unique irreducible $\ih(n-1)$-module of dimension $(n-2)$.
\end{Lemma}
\noin\textbf{Proof.} The existence part is provided by the module
$S$ above. The uniqueness part is more difficult and is treated in
Proposition $17$ page $145$ of \cite{AR}.\\

\noin This ends the proof of Theorem $7$. We note that this theorem remains true for $n=4$. \\
The next theorem describes the irreducible
$\dbw$-dimensional invariant subspace of $\V^{(n)}$ when $l=-r^3$.
\begin{theo}
Let $n$ be an integer with $n\geq 4$ and assume $\ih(n)$ is
semisimple. Suppose $l=-r^3$.\\
When $n=4$, the irreducible $3$-dimensional invariant subspace of
$\V^{(4)}$ is spanned by the vectors $u_1$, $u_2$, $u_3$ of Theorem
$4$.\\
When $n\geq 5$, the irreducible $\dbw$-dimensional invariant
subspace of $\V^{(n)}$ is built inductively from the irreducible
$\frac{(n-2)(n-3)}{2}$-dimensional invariant subspace of
$\V^{(n-1)}$ by adding the $(n-2)$ linearly independent vectors:
$$V_k^{(n)}=w_{k+1,n}-r\,w_{k,n}+r^{n-k}\,w_{k,k+1},\;\; k=1,\dots,n-2$$
\end{theo}
\noin\textbf{Proof.} First, we show a lemma.
\begin{Lemma}
Let $n\geq 5$. The vectors $V_1^{(n)}$,$\dots$,$V_{n-2}^{(n)}$
belong to $K(n)$. Moreover, the action of the $g_i$'s on these
vectors is as follows:
\begin{eqnarray*}
g_{k-1}.V_k^{(n)}&=&
V_{k-1}^{(n)}+r\,V_k^{(n)}-r^{n-k-1}\,V_{k-1}^{(k+1)}\\
g_k.V_k^{(n)}&=&-\unsurr\,V_k^{(n)}\\
g_{k+1}.V_k^{(n)}&=&V_{k+1}^{(n)}+r\,V_k^{(n)}-r^{n-k-1}\,V_k^{(k+2)}\;\;\text{when $k<n-2$}\\
g_{n-1}.V_{n-2}^{(n)}&=&-\unsurr\,V_{n-2}^{(n)}\\
g_{n-1}.V_k^{(n)}&=&V_k^{(n-1)}-m\,V_k^{(n)}\;\;\;\;\text{when $k<n-2$}\\
g_i.V_k^{(n)}&=&r\,V_k^{(n)}\;\;\;\,\;\;\qquad\qquad\text{when
$i\not\in\lb k-1,k,k+1,n-1\rb$}
\end{eqnarray*}

%and span an irreducible $\ih(n-1)$-module of
%dimension $(n-2)$.
\end{Lemma}
\noin \textbf{Proof.} The fact that these vectors belong to $K(n)$ is Claim $3$ page $120$ of \cite{AR}.
The equalities that follow are obtained by straightforward computations.\\\\
%Let's compute the action of $g_1$,$\dots$,$g_{n-2}$ on the vectors
%$V_k^{(n)},\;k=1,\dots,n-2$. Fix $k$ with $1\leq k\leq n-2$.
%Straightforward computations show that:
%In the other cases, the action by $g_i$ on $V_k^{(n)}$ is a
%multiplication by $r$.
%\begin{Lemma}
%Let $n$ be an integer with $n\geq 5$. Assume that $\ih(n)$ is
%semisimple. Suppose $l=-r^3$ and $r^{2(n-1)}=-1$. In $K(n)$ there
%exists a unique irreducible $\ih(n-1)$-module of dimension $(n-2)$.
%\end{Lemma}
%\noin \textbf{Proof.} Let's show the existence of such a module.
%Since $r^{2n}\neq -1$, $K(n)$ is irreducible and $\dbw$-dimensional.
%We have $0\subset K(n)\cap\V^{(n-1)}\subseteq K(n-1)$. Moreover, we
%know that $K(n-1)$ is a direct sum of a one-dimensional invariant
%subspace and of an irreducible $\frac{(n-2)(n-3)}{2}$-dimensional
%invariant subspace of $\V^{(n-1)}$. Further, by Lemma $13$ page
%$139$ of \cite{AR}, there does not exist any one-dimensional
%$\ih(n-1)$-module inside $K(n)$. Then $K(n)\cap\V^{(n-1)}$ is
%irreducible, $\frac{(n-2)(n-3)}{2}$-dimensional. Let $S$ be an
%$\ih(n-1)$-module of dimension $(n-2)$ that is a summand of
%$K(n)\cap\V^{(n-1)}$ in $K(n)$. Since by Lemma $13$ page $139$ of
%\cite{AR}, there does not exist any one-dimensional
%$\ih(n-1)$-module in $K(n)$, Corollary $1$ forces
%that $S$ is irreducible.\\\\
\indent Let's prove Theorem $8$. Suppose $n\geq 5$. We distinguish
between several cases.
\begin{enumerate}
\item If $r^{2(n-1)}\neq -1$ and $r^{2n}\neq -1$.\\\\
Then $K(n-1)$ is irreducible and has dimension
$\frac{(n-2)(n-3)}{2}$. Also $K(n)$ is the irreducible
$\dbw$-dimensional invariant subspace of $\V^{(n)}$.  The
irreducibility of $K(n-1)$ and the fact that $\dbw>n-1$ when $n\geq
5$ imply that $K(n-1)\subseteq K(n)$. Since we notice that
$k(n)=k(n-1)+(n-2)$, $K(n)$ is a direct sum of $K(n-1)$ and of an
$(n-2)$-dimensional vector space spanned by the vectors $V_1^{(n)}$,
$V_2^{(n)}$, $\dots$, $V_{n-2}^{(n)}$.
\item If $r^{2(n-1)}=-1$.\\\\
Then $r^{2n}\neq -1$. So $K(n)$ is the irreducible
$\dbw$-dimensional invariant subspace of $\V^{(n)}$. We have
$0\subset K(n)\cap\V^{(n-1)}\subseteq K(n-1)$. Moreover, we know
that $K(n-1)$ is a direct sum of a one-dimensional invariant
subspace and of an irreducible $\frac{(n-2)(n-3)}{2}$-dimensional
invariant subspace of $\V^{(n-1)}$. Further, since by Lemma $13$
page $139$ of \cite{AR}, there does not exist any one-dimensional
$\ih(n-1)$-module inside $K(n)$, we see that $K(n)\cap\V^{(n-1)}$ is
irreducible, $\frac{(n-2)(n-3)}{2}$-dimensional. So,
$K(n)\cap\V^{(n-1)}$ is the irreducible
$\frac{(n-2)(n-3)}{2}$-dimensional invariant subspace of
$\V^{(n-1)}$ and $K(n)$ is a direct sum of $K(n)\cap\V^{(n-1)}$ and
of an $(n-2)$-dimensional vector space spanned by the vectors
$V_1^{(n)}$, $V_2^{(n)}$, $\dots$, $V_{n-2}^{(n)}$.
\item If $r^{2n}=-1$.\\\\
Let $\W$ be the irreducible $\dbw$-dimensional invariant subspace of
$\V^{(n)}$. We have seen along the proof of Theorem $7$ that $\W$
contains $K(n-1)$, the irreducible
$\frac{(n-2)(n-3)}{2}$-dimensional invariant subspace of
$\V^{(n-1)}$. Consider the $F$-vector space
$$\mathcal{S}=K(n-1)\oplus\,\text{Span}_F(V_1^{(n)},\dots,V_{n-2}^{(n)})$$
This vector space has dimension $\frac{(n-1)(n-2)}{2}$ over $F$. We
will show that it is stable under the action by the $g_i$'s and that
$\W=\mathcal{S}$. When $r^{2n}=-1$, we have $r^2\neq -1,(r^2)^2\neq
-1,\dots,r^{2(n-1)}\neq -1$. We notice that $K(3)$ is spanned over
$F$ by $V_1^{(3)}$ (see Theorem $3$) and a use of the computer
program of appendix $A$ of \cite{THE} shows that the vectors
$V_1^{(3)}$, $V_1^{(4)}$ and $V_2^{(4)}$ belong to $K(4)$. In
particular, we have $K(3)\subset K(4)$. These remarks and point $1.$
of the ongoing proof imply that
$$(\divideontimes)\left\lb\begin{array}{l}
K(3)\subset K(4)\subset\dots\subset K(n-1),\\
K(s)=\,<V_i^{(j)}: 3\leq j\leq s,\,1\leq i\leq j-2>_{F}
\end{array}\right.$$
Further, we have:\begin{center}
$\begin{array}{ccc}g_{n-1}.V_k^{(j)}&=&r\,V_k^{(j)}\\
g_{n-1}.V_k^{(n-1)}&=&V_k^{(n)}\end{array}$
$\begin{array}{l}\forall\,3\leq j\leq n-2,\;\forall\,1\leq k\leq
j-2\\\forall\,1\leq k\leq n-3\end{array}$\end{center} It follows
that
$$g_{n-1}.\,K(n-1)\subseteq\mathcal{S}$$
This inclusion, point $(\divideontimes)$ and the equalities of Lemma
$4$ imply that $\mathcal{S}$ is stable under the action by the
$g_i$'s. Since $\mathcal{S}$ is contained in $K(n)$, it follows that
$\mathcal{S}$ is an $\ih(n)$-module. Moreover, by choice of $l$ and
$r$, we see that $\mathcal{S}$ must be irreducible. Hence
$\mathcal{S}$ is the irreducible $\dbw$-dimensional invariant
subspace of $\V^{(n)}$. This ends the proof of Theorem $8$.
\end{enumerate}
\section{Proof of the Main Theorem}
We will work on the representation $\n^{(n)}$ of this paper instead
of the original representation of Krammer of \cite{KR}. This is
allowed by the following lemma.
\begin{Lemma}
It suffices to check the main theorem on the representation
$\n^{(n)}$ where $q$ has been replaced by $\unsur{r^2}$ and $t$ has
been replaced with $\frac{r^3}{l}$.
\end{Lemma}
\noindent\textbf{Proof.} We recall from $\S\,2$ that the
representation of this paper is equivalent to the Lawrence--Krammer
representation of the BMW algebra of \cite{CGW} where our $r$ is the
$\unsurr$ of \cite{CGW}. Further, up to some rescaling of the
generators, the representation of \cite{CGW} is equivalent, as a
representation of the braid group on $n$ strands, to the
representation of the Artin group of type $A_{n-1}$ of \cite{CW}.
The parameters $t$ and $r$ of \cite{CW} are related to the
parameters $l$ and $r$ of \cite{CGW} by $l=\unsur{tr^3}$. The
representation of \cite{CW} is itself equivalent to the original
representation of Krammer of \cite{KR} with parameters $t$ and $q$.
The link between the parameter $q$ of \cite{KR} and the parameter
$r$ of \cite{CW} is given by $q=r^2$. \\%Hence it suffices to check

%main theorem on our representation
%where $q$ has been replaced by $\unsur{r^2}$ and $t$ has been
%replaced with $\frac{r^3}{l}$.

\textsc{Proof of the Main Theorem.\\\\} \indent We first deal with
the uniqueness part when we exclude the case when $l=-r^3$ and
$r^{2n}=-1$. Then, under the assumption that $\ih(n)$ is semisimple,
the values $\unsur{r^{2n-3}}$, $\unsur{r^{n-3}}$,
$-\unsur{r^{n-3}}$, $r$, $-r^3$ are all distinct. When $l=r$ or
$l=-r^3$, we have seen that $K(n)$ is irreducible. As any proper
invariant subspace of $\V^{(n)}$ must be contained in $K(n)$ (this
is Remark $2$), the $B(A_{n-1})$-module $K(n)$ is then the unique
proper invariant subspace of $\V^{(n)}$. Next, suppose that
$l\in\lb\jca,-\jca,\unsur{r^{2n-3}}\rb$ and let $\W$ be an
irreducible invariant subspace of $\V^{(n)}$. If
$\text{dim}(\W)\geq\dbw$, then $\text{dim}(\W)>n-1$ as soon as
$n\geq 5$. This implies that $\W\cap\V^{(n-1)}\neq \lb 0\rb$, hence
$\n^{(n-1)}$ is reducible. Also, we have $\text{dim}(\W)> 2n-3$ as
soon as $n\geq 6$. This implies that $\W\cap\V^{(n-2)}\neq\lb 0\rb$,
hence $\n^{(n-2)}$ is reducible. Using the reducibility theorem of
$\S 2$ and our usual restriction on the parameter $r$, we then get
$l\in\lb r,-r^3\rb$, a contradiction. By the same proof when $n=8$,
we cannot have $\text{dim}(\W)=14$ since $14>13$. So assuming $n\geq
6$ and using Corollary $1$ and Theorem $5$, we must hence have
$\text{dim}(\W)\in\lb1,n-1\rb$. Now the uniqueness follows from
Theorem $3$ and Theorem $4$. When $n=5$, the irreducible
representations of $\ih(5)$ have degrees $1$, $4$, $5$ and $6$. If
$\text{dim}(\W)=6$, applying the reducibility theorem on $\n^{(4)}$,
we get $l\in\lb r,-r^3,\unsurr,-\unsurr,\unsur{r^{5}}\rb$. Since we
assumed $l\in\lb\unsur{r^7},\unsur{r^2},-\unsur{r^2}\rb$, we then
have $l=\unsur{r^7}$ and $l\in\lb\unsurr,-\unsurr\rb$. This means
that there exists a (unique) one-dimensional invariant subspace in
$\V^{(5)}$ and an irreducible $3$-dimensional invariant subspace in
$\V^{(4)}$, which is also the unique proper invariant subspace of
$\V^{(4)}$. We get $k(5)=7$ and so $\di(K(5)\cap\V^{(4)})\geq 3$.
But $\V^{(4)}$ has a unique proper invariant subspace which has
dimension $3$, so we get $K(5)\cap\V^{(4)}=K(4)$. In particular, we
derive $K(4)\subseteq K(5)$. Further, by Theorem $4$,
$$v_2^{(4)}=\bigg(\unsurr-\unsur{l}\bigg)\,w_{23}+\big(w_{24}-\unsurr\,w_{34}\big)+\varepsilon_{l}\,r\,\big(w_{12}-\unsurr\,w_{13}\big)$$
belongs to $K(4)$. We show that the action of $C_{35}$ on
$v_2^{(4)}$ is nonzero. The action of the elements $C_{ij}$'s is
summarized in appendix $C$ of \cite{THE}. We have
\begin{eqnarray*}
C_{35}.\,w_{23}&=&\unsurr\,w_{35}\\
C_{35}.\,w_{24}&=&\bigg(\unsurr-r\bigg)\bigg(\unsur{l}-\unsurr\bigg)\,w_{35}\\
C_{35}.\,w_{34}&=&\unsur{l}\,w_{35}\\
C_{35}.\,w_{13}&=&\unsur{r^2}\,w_{35}
\end{eqnarray*}
It follows that $$\begin{array}{ccccccccc} \text{If}& l=\unsurr,
&\text{then}& C_{35}.\,v_2^{(4)}&=&(-r^2-\unsur{r^2})\,w_{35}&\neq
0&
\text{since} &r^4\neq-1.\\
\text{If} &l=-\unsurr,& \text{then}&
C_{35}.\,v_2^{(4)}&=&(r+\unsurr)^2\,w_{35}&\neq
0&\text{since}&r^2\neq -1.\end{array}$$ These computations show that
$v_2^{(4)}\not\in K(5)$, which contradicts $K(4)\subseteq K(5)$. We
conclude that it is impossible to have $\di(\W)=6$. Thus, an
invariant subspace is unique in the case when $n=5$ as well. This
ends the uniqueness part in point number $1$.
\\%and $\di
%(K(5)\cap\V^{(4)})=3$. Let $S$ be an
%$\ih(4)$-module that is a summand of $K(5)\cap\V^{(4)}$ in $K(5)$.
%Hence $\di (S)=4$. This dimension is incompatible with the values
%for $l$ and $r$. This ends the proof of uniqueness.\\
\indent Let's study the case when $l=-r^3$ and $r^{2n}=-1$. As
$K(n)$ is a direct sum of the unique one-dimensional invariant
subspace and of the unique irreducible $\dbw$-dimensional invariant
subspace of $\V^{(n)}$, there are indeed exactly three proper
invariant subspaces in $\V^{(n)}$. It remains to show that when
$n\geq 5$, the unique irreducible $\dbw$-dimensional invariant
subspace of $\V^{(n)}$ is isomorphic to the Specht module
$S^{(n-2,1,1)}$. This is the object of Proposition $18$ page $149$
of \cite{AR}. We recall that as part of the proof of Theorem $3$,
the Specht module $S^{(1,1,1)}$ occurs in $\V^{(3)}$ when $l=-r^3$
and as part of the proof of Theorem $4$, the Specht module
$S^{(2,1,1)}$ occurs in $\V^{(4)}$ when $l=-r^3$. Point number $2)$
is thus
entirely proven.\\
\indent Let's go back to point number $1)$. The unique proper
invariant subspace of $\V^{(n)}$ is $K(n)$. We call it $\W$. When
$l=\unsur{r^{2n-3}}$, $\W$ is isomorphic to $S^{(n)}$ by the proof
of Theorem $3$. Suppose next $l\in\lb\jca,-\jca\rb$. Then $\W$ is
isomorphic to $S^{(n-1,1)}$ by the proof of Theorem $4$. Suppose now
$l=r$. By Theorem $5$, $\W$ is an irreducible $\cil$-dimensional
$\ih(n)$-module. By Theorem $2$, such an irreducible $\ih(n)$-module
must be isomorphic to $S^{(n-2,2)}$ or its conjugate Specht module
$S^{(2,2,1^{n-4})}$, except when $n=7$ in which case it can also be
isomorphic to $S^{(4,3)}$ or its conjugate Specht module
$S^{(2,2,2,1)}$, both of dimension $14$. Thus, there are two things
to show. One of them is that for $n\geq 5$, $\W$ cannot be
isomorphic to $S^{(2,2,1^{n-4})}$. The other one is that $\W$ cannot
be isomorphic to $S^{(4,3)}$ or to $S^{(2,2,2,1)}$ when $n=7$. Let's
start with the first point. We proceed by induction on $n\geq 5$ to
show that if $\W$ is an irreducible $\cil$-dimensional invariant
subspace of $\V^{(n)}$, it is impossible to have $\W\simeq
S^{(2,2,1^{n-4})}$. We show that this is true when $n=5$. This is
Result $1$ page $21$ of \cite{CLW}. When $n\geq 6$, we use the
branching rule. If $\W\simeq S^{(2,2,1^{n-4})}$, then
$$\W\da_{\ih(n-1)}\simeq S^{(2,1^{n-3})}\oplus S^{(2,2,1^{n-5})}$$
We have $0\subset\W\cap\V^{(n-1)}\subset\W\da_{\ih(n-1)}$, so we get % Further,
%we have
%$$\di(\W\cap\V^{(n-1)})>\cil-(n-1)=\frac{n(n-5)}{2}+1$$
%So,
%$$\di(\W\cap\V^{(n-1)})\geq\frac{n(n-5)}{2}+2=\frac{(n-1)(n-4)}{2}$$
$\W\cap\V^{(n-1)}\simeq S^{(2,2,1^{n-5})}$. This is in contradiction
with our induction hypothesis and we are done with the first point.
Let's deal with the second point. First, if $\W$ is isomorphic to
$S^{(2,2,2,1)}$, then by the branching rule,
$$\W\da_{\ih(6)}\simeq S^{(2,2,2)}\oplus S^{(2,2,1,1)}$$
Since by the proof of Theorem $4$, the Specht module $S^{(2,2,2)}$
cannot occur in $\V^{(6)}$, we must have $\W\cap\V^{(6)}\simeq
S^{(2,2,1,1)}$. Another application of the branching rule then
yields
$$\W\cap\V^{(6)}\da_{\ih(5)}\simeq S^{(2,1,1,1)}\oplus S^{(2,2,1)}$$
The intersection $\W\cap\V^{(5)}$ is nonzero. Moreover,
$S^{(2,1,1,1)}$ does not occur in $\V^{(5)}$ by the proof of Theorem
$4$. Also $S^{(2,2,1)}$ does not occur in $\V^{(5)}$ (this is Result
$1$ page $21$ of \cite{CLW}). We thus get a contradiction. If now
$\W$ is isomorphic to $S^{(4,3)}$, applying the branching rule twice yields %we must have
%$\W\cap\V^{(6)}\simeq S^{(4,2)}$ since $S^{(3,3)}$ cannot occur in
%$\V^{(6)}$ by the proof of Theorem $4$.
\begin{equation}\W\da_{\ih(5)}\simeq 2\,S^{(3,2)}\oplus S^{(4,1)}\end{equation} Let $w_1$,
$w_2$, $w_3$, $w_4$, $w_5$ be linearly independent vectors of $\W$
such that the actions by $g_1$, $g_2$, $g_3$ and $g_4$ on these
vectors is given by the matrices $P_i$'s of Fact $1$ page $19$ of
\cite{CLW}. Notice the presence of a seventh node does not change
the computations of \cite{CLW}. Indeed, since
$g_1.w_4=-\unsurr\,w_4$ and $g_3.\,w_4=-\unsurr\,w_4$, we see that
the vector $w_4$ belongs to $\V^{(4)}$. Then the relations
$w_5=g_4.\,w_4$, $w_1=g_2.\,w_4-r\,w_4$, $w_2=g_4.\,w_1$,
$w_3=g_3.\,w_2-r\,w_2$, show that the presence of a seventh node can
just be forgotten. Then, up to a multiplication by a scalar, the
$w_i$'s are uniquely determined as in \cite{CLW}. This contradicts
the multiplicity of $S^{(3,2)}$ in $(1)$. So we are done with the
case $l=r$. It remains to deal with the reducibility case $l=-r^3$
when we assume $r^{2n}\neq -1$. In this case $\W$ is an irreducible
$\dbw$-dimensional $\ih(n)$-module. A use of Theorem $2$ shows that
$\W$ must be isomorphic to $S^{(n-2,1,1)}$ or to $S^{(3,1^{n-3})}$.
We show that $\W$ is isomorphic to $S^{(n-2,1,1)}$. This is already
true for $n=3$ by the proof of Theorem $3$ and for $n=4$ by the
proof of Theorem $4$. When $n=5$, there is nothing to prove as
$S^{(3,1,1)}$ is self-conjugate. As for $n\geq 6$, we proceed by
induction on $n$ to show that $S^{(3,1^{n-3})}$ cannot occur in
$\V^{(n)}$. First we show that $S^{(3,1,1,1)}$ cannot occur in
$\V^{(6)}$. Suppose $K(6)\simeq S^{(3,1,1,1)}$. Then by the
branching rule,
$$K(6)\da_{\ih(4)}\,\simeq\, 2\,S^{(2,1,1)}\oplus S^{(3,1)}\oplus
S^{(1,1,1,1)}$$
%$$K(6)\da_{\ih(5)}\simeq S^{(3,1,1)}\oplus
%S^{(2,1,1,1)}$$ By inspection on the dimensions, we thus get
%$K(6)\cap\V^{(5)}\simeq S^{(3,1,1)}$. %In particular,
%$K(6)\cap\V^{(5)}$ is the irreducible $6$-dimensional invariant
%subspace of $\V^{(5)}$.
%In particular, it will be
%useful to note that $K(6)\cap\V^{(5)}$ is the irreducible
%$6$-dimensional invariant subspace of $\V^{(5)}$.
%By the branching rule again, we have
%$$K(6)\cap\V^{(5)}\da_{\ih(4)}\simeq S^{(2,1,1)}\oplus S^{(3,1)}$$
%From there, we derive
%$K(6)\cap\V^{(4)}=\text{Span}_{F}(u_1,u_2,u_3)$. Now
%$K(6)\cap\V^{(4)}$ is the irreducible $3$-dimensional invariant
%subspace of $\V^{(4)}$ and $K(6)\cap\V^{(5)}$ is the irreducible
%$6$-dimensional invariant subspace of $\V^{(5)}$, so that by Theorem
%$8$ we get
%$$K(6)\cap\V^{(5)}=\text{Span}_{F}(u_1,u_2,u_3)\oplus\text{Span}_{F}(V_1^{(5)},V_2^{(5)},V_3^{(5)})$$
%Let $\mathcal{S}$ be an $\ih(5)$-summand of $K(6)\cap\V^{(5)}$ in
%$K(6)$. This summand $\mathcal{S}$ is isomorphic to $S^{(2,1,1,1)}$.
%restriction to $\ih(4)$ is isomorphic to $S^{(1,1,1,1)}\oplus
%S^{(2,1,1)}$.
In particular, there exists a vector $w$ of $K(6)$ such that
$$g_i.\,w=-\unsurr\,w\;\;\text{for each $i\in\lb 1,2,3\rb$}$$ These
relations imply that in $w$ there are no terms in $w_{1j}$ or
$w_{2j}$ or $w_{3j}$ or $w_{4j}$ or $w_{jk}$ for any $j\geq 5$ and
any $k\geq j+1$. In other words, $w$ belongs to $\V^{(4)}$. But in
$\V^{(4)}$ there does not exist any one-dimensional invariant
subspace isomorphic to $S^{(1^4)}$, hence a contradiction. This
finishes the base case $n=6$. Let $n\geq 7$ and suppose that
$S^{(3,1^{n-4})}$ cannot occur inside $\V^{(n-1)}$. Let's show that
$S^{(3,1^{n-3})}$ cannot occur inside $\V^{(n)}$. If $\W$ is an
irreducible $\dbw$-dimensional invariant subspace of $\V^{(n)}$ that
is isomorphic to $S^{(3,1^{n-3})}$, its restriction to $\ih(n-1)$ is
isomorphic to $S^{(3,1^{n-4})}\oplus S^{(2,1^{n-3})}$. Then
$\W\cap\V^{(n-1)}$ is isomorphic to $S^{(3,1^{n-4})}$. This is
impossible with our induction hypothesis. This ends the proof of the
Main Theorem.


\begin{thebibliography}{lab}
\bibitem{BI} S.J. Bigelow, The Burau representation is not faithful
for $n=5$, \textit{Geometry and Topology} \textbf{3} $(1999)$,
$397-404$
\bibitem{BIG} S.J. Bigelow, Braid groups are linear, \textit{J.
Amer. Math. Soc.}, \textbf{14} $(2001)$, $471-486$
\bibitem{BW} J.S. Birman and
H. Wenzl, Braids, link polynomials and a new algebra \textit{Trans.
Amer. Math. Soc.} \textbf{313} $(1989)$, no. $1$, $249-273$
\bibitem{CW} A.M. Cohen and D.B. Wales, Linearity of Artin groups of
finite type, \textit{Isr. J. Math.}, \textbf{131} $(2002)$,
$101-123$
\bibitem{CGW}
A.M. Cohen, D.A.H. Gijsbers and D.B. Wales, BMW algebras of simply
laced type, \textit{J. Algebra}, \textbf{285} $(2005)$, no.$2$,
$439-450$
\bibitem{DI}
F. Digne, On the linearity of Artin braid groups, \textit{J.
Algebra}, \textbf{268}, No.1, $(2003)$, $39-57$
\bibitem{JAM} G.D. James, On the minimal dimensions of
irreducible representations of symmetric groups, \textit{Math. Proc.
Camb. Phil. Soc.} \textbf{94} $(1983)$, $417-424$
\bibitem{KR} D. Krammer, Braid groups are linear, \textit{Ann. of Math.}
\textbf{155} $(2002)$, $131-156$
\bibitem{RUT} R. Lawrence, Homological representations of the Hecke
algebras, \textit{Comm. Math. Physics}, \textbf{135}, $(1990)$,
$141-191$
%\bibitem{DIM} C. Levaillant, Dimensions of the invariant subspaces
%of the Lawrence--Krammer representation, in preparation
\bibitem{THE} C. Levaillant, Irreducibility of the Lawrence--Krammer
representation of the BMW algebra of type $A_{n-1}$, Ph.D. thesis
California Institute of Technology $(2008)$
http://thesis.library.caltech.edu/2255/1/thesis.pdf
\bibitem{AR} C. Levaillant, Irreducibility of the Lawrence--Krammer
representation of the BMW algebra of type $A_{n-1}$, arXiv:0901.3908
[arXiv version of the Ph.D. thesis]
\bibitem{LEV} C. Levaillant, Irreducibility of the Lawrence--Krammer
representation of the BMW algebra of type $A_{n-1}$, \textit{C.R.
Acad. Sci. Paris, Ser. I} $347$ $(2009)$, $15-20$
\bibitem{CLW} C.I. Levaillant and D.B. Wales, Parameters for which the
Lawrence--Krammer representation is reducible, arXiv:0901.3856.
[Shorter version available in \textit{J. Algebra}, $323$ $(2010)$,
$1966-1982$]
\bibitem{dim} C. Levaillant, Dimensions of the invariant subspaces
of the Lawrence--Krammer representation when it is reducible, in
preparation
\bibitem{LP} D.D. Long and M. Paton, The Burau representation is not
faithful for $n\geq 6$, Topology \textbf{32} (1993) no.2, $439-447$
\bibitem{MAR} I. Marin, Sur les repr\'esentations de Krammer
g\'en\'eriques, \textit{Ann. Inst. Fourier}, \textbf{57}, no.6,
$(2007)$, $1883-1925$
\bibitem{MAT} A. Mathas, Iwahori-Hecke algebras and Schur algebras of the
symmetric group, \textit{University Lecture Series}, Volume
\textbf{15}.
\bibitem{MOO} J.A. Moody, The Burau representation of the braid group $B_n$ is
unfaithful for large $n$, \textit{Bull. Amer. Math. Soc.}
\textbf{25} $(1991)$, no. $2$, $379-384$
\bibitem{MW} H.R. Morton and A.J. Wasserman, A basis for the
Birman-Wenzl algebra, preprint $1989$
\bibitem{MUR} J. Murakami, The Kauffman polynomial of links and
representation theory, \textit{Osaka J. Math.} \textbf{24} $(1987)$,
$745-758$
\bibitem{RUI} H. Rui and M. Si, Blocks of Birman-Murakami-Wenzl algebras, \textit{Int. Math. Res. Not.} \textbf{rnq083}
$(2010)$
\bibitem{Z} M. Zinno, On Krammer's representation of the
braid group, \textit{Math. Ann.} $321$ $(2001)$, $197-211$
\end{thebibliography}
\end{document}